\documentclass[11pt]{amsart}
\usepackage[margin=1.0in]{geometry}
\usepackage{amsmath,amssymb,amsfonts,graphicx,color,fancyhdr,psfrag,comment,enumerate}
\usepackage[hyperpageref]{backref}
\usepackage[colorlinks=true, pdfstartview=FitV, linkcolor=blue, citecolor=blue, urlcolor=blue]{hyperref}

\usepackage[nodisplayskipstretch]{setspace}
\usepackage[utf8]{inputenc}

\usepackage{epstopdf}
\usepackage{mathrsfs}
\usepackage{epsfig}
\usepackage{hyperref}
\usepackage{ifthen}
\usepackage{float}
\usepackage{enumerate}
\usepackage{mathtools}
\usepackage[bottom]{footmisc}
\usepackage{tikz}
\usetikzlibrary{matrix,shapes,arrows,positioning,chains}
\usepackage{amssymb,amsmath,amsfonts}
\usepackage{bbm}
\usepackage[dvipsnames]{xcolor}
\usepackage[normalem]{ulem}

\numberwithin{equation}{section}

\newcommand{\R}{\mathbb R}
\newcommand{\N}{\mathbb N}

\newcommand{\mf}[1]{\mathfrak{#1}}

\newcommand{\C}{\mathbb C}

\newtheorem{theorem}{Theorem}[section]
\newtheorem{prop}[theorem]{Proposition}
\newtheorem{lemma}[theorem]{Lemma}

\newtheorem{remark}[theorem]{Remark}

\newcommand{\Ha}{\mathbb{H}^n}

\usepackage{tikz}
\usepackage{pgfplots}
\pgfplotsset{compat=1.13,
tick label style={color=white},
  label style={font=\small},
  legend style={font=\small}}

\usepackage{mathtools}

\allowdisplaybreaks

\newcommand*{\Resize}[2]{\resizebox{#1}{!}{$#2$}}%

\title[Bilinear Kor\'anyi spherical maximal function]{Bilinear Kor\'anyi spherical maximal function on the Heisenberg group}

\author[Ghosh]{Abhishek Ghosh}
\address[Abhishek Ghosh]{Department of Mathematics, Indian Institute of Technology Madras, Chennai, 600036, India.}

\email{abhi@iitm.ac.in, abhi170791@gmail.com}

\author[Singh]{Rajesh K. Singh}
\address[Rajesh K. Singh]{Department of Mathematics \& Statistics, Gitam Institute of Science, GITAM University, Visakhapatnam- 530045, A.P., India}
\email{rsingh4@gitam.edu}



\begin{document}

\begin{abstract}
We introduce the  bilinear Kor\'anyi spherical means on the Heisenberg group $\mathbb{H}^n,$ and derive $L^{p_1}(\mathbb{H}^n) \times L^{p_2}(\mathbb{H}^n) \to L^{p}(\mathbb{H}^n)$ estimates for the associated bilinear Kor\'anyi spherical maximal function on $\mathbb{H}^n.$ We also address multilinear extensions of our results. 
\end{abstract}
	
\maketitle	


\section{Introduction and Preliminaries}
The study of spherical maximal functions is a classical problem in Harmonic analysis and goes back to Stein's spherical maximal theorem, \cite{SteinSph}, which asserts that the spherical maximal operator 
$$S_{\R^n}f(x):=\sup_{r>0}\bigg|\int_{S^{n-1}} f(x-ry)\, d\sigma_{n-1}(y)\bigg|,$$ is bounded from $L^p(\R^n)$ to itself if and only if $p>\frac{n}{n-1}, n\geq 3.$ Here, $d \sigma_{n-1}$ is the rotationally invariant measure on the sphere $S^{n-1}$, normalized such that $\sigma_{n-1}(S^{n-1})=1$ . Later, in the celebrated work \cite{Bourgain}, Bourgain proved the result in dimension two, see also \cite{MSS}. In recent years, considerable attention has been devoted to extending Stein's theorem beyond the Euclidean setting. A pioneering work of Cowling \cite{Cowling} established the analogue of Stein's spherical maximal theorem on the Heisenberg group. Before stating his result, we briefly recall some preliminaries on the Heisenberg group.

Let $ \mathbb{H}^n:= \mathbb{C}^n \times \mathbb{R} $ denotes the $d:=(2n + 1)$-dimensional Heisenberg group, equipped with the group law defined for $ x, y \in \mathbb{H}^n $ by  
\[
x \cdot y := \left( \underline{x} + \underline{y},~ \overline{x} + \overline{y} + \textstyle\frac{1}{2}  \Im(  \underline{x} \cdot \bar{\underline{y}}) \right). 
\]  
 We denote the non-central component of $ x \in \mathbb{H}^n $ by $ \underline{x} \in \mathbb{C}^n $ and the central component by $ \overline{x} \in \mathbb{R},$ sometimes we shall also use $x=(z, t),$ with $z\in \C^n, t\in \R.$ The Haar measure on $\mathbb{H}^n$ coincides with the Lebesgue measure $dx=d\underline{x}\,d\overline{x}$, and $\Ha$ is equipped with a left-invariant metric induced by
the Kor\'anyi norm  defined by
 \begin{equation*}
     |x|:= \left( \|\underline{x} \|^4 + |\overline{x}|^2 \right)^{\frac{1}{4}}, \quad\quad x=(\underline{x}, \overline{x}) \in \mathbb{H}^n,
 \end{equation*}
where $\| \underline{x}  \|$ denotes the Euclidean norm of $\underline{x} \in \mathbb{C}^{n}$. Notably, $(\Ha, |\cdot|)$ is a space of homogeneous type. The Kor\'anyi norm is homogeneous of degree one with respect to the anisotropic family of dilations given by $\delta_{r}(\underline{x}, \overline{x}):=(r\underline{x}, r^2\overline{x}),$ that is, $|\delta_{r}(\underline{x}, \overline{x})|=r|(\underline{x}, \overline{x})|$,\, $\forall \,\, 0<r<\infty , \, \forall \, (\underline{x}, \overline{x}) \in \Ha$. For a function $f$  on $\Ha,$ we denote its anisotropic dilation, $\delta_r f $, by $ \delta_r f(\underline{x},\overline{x})=f(\delta_r(\underline{x},\overline{x})),\ \, (\underline{x}, \overline{x})\in \Ha.$ For a compactly supported Borel measure $\mu$ on $\Ha,$ its dilates, $\delta_r \mu$, are defined as $\langle \delta_r \mu, \phi \rangle :=\langle \mu, \delta_r \phi\rangle,\ \forall \, \phi\in C^{\infty}_{c}(\Ha).$ Let $\Sigma_{K}:=\{x\in \Ha: 
|x|=1\}$ be the Kor\'anyi sphere and $\mu$ be the unique Radon measure on $\Sigma_{K}$ satisfying the polar decomposition 
\[\int_{\Ha} f(x)\, dx=\int_{0}^{\infty}\int_{\Sigma_{K}}f(\delta_{r}x) r^{Q-1}\, d\mu(\omega)\, dr,\]
for all $f\in L^1(\Ha).$ Here, $Q=2n+2$ is the homogeneous dimension of $\Ha.$ Set $\mu_r = \delta_r \mu,$ then the Kor\'anyi spherical means are defined as $$ A^{K}_rf(x) = f\ast \mu_r(x) = \int_{|y|=1} f(x.\delta_r y^{-1}) d\mu(y),$$
and the Kor\'anyi maximal function is $$S_{\Ha}f(x) = \sup_{r>0} |A^{K}_rf(x)|, \quad x\in \Ha.$$
In \cite{Cowling}, M. Cowling proved the following sharp result.
\begin{theorem} The Kor\'anyi spherical maximal function $ S_{\Ha} $ is bounded on $ L^p(\Ha) $ for all $ p > \frac{2n+1}{2n}, n\geq 2.$
\label{Koranyi:Thm}
\end{theorem}
The primary purpose of this article is to introduce the bilinear Kor\'anyi spherical means, and study the Lebesgue space estimates of the associated bilinear Kor\'anyi maximal operator. We now delve into defining the main objects of this article. 
\medskip

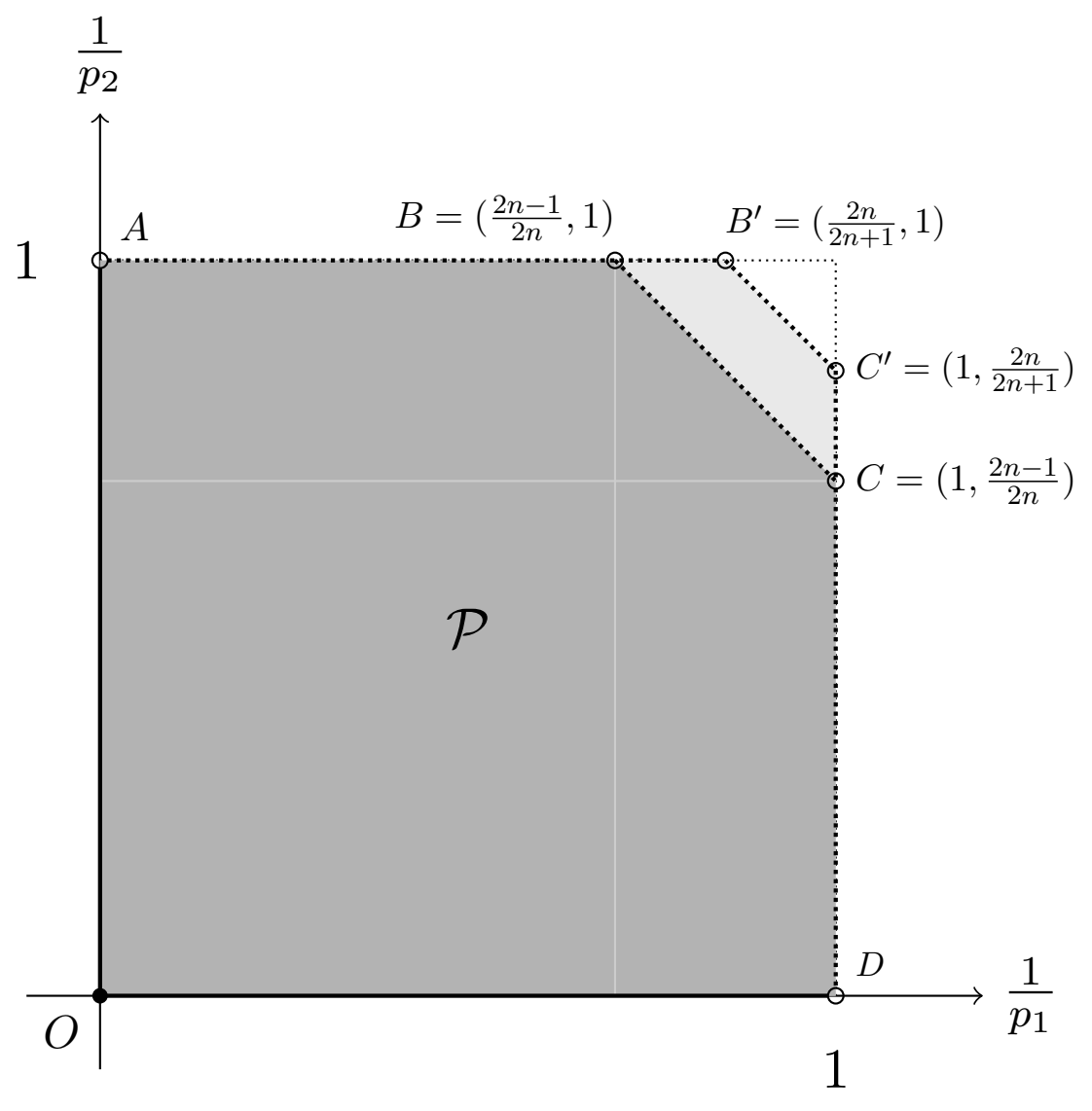
\begin{figure}[ht]
\centering
\begin{tikzpicture}[scale=5] 
    \draw[->] (-0.1,0) -- (1.2,0) node[right] {$\frac{1}{p_{1}}$};
    \draw[->] (0,-0.1) -- (0,1.2) node[above] {$\frac{1}{p_{2}}$};
    
    \draw[densely dotted] (0,1) -- (1,1) -- (1,0);
    \node at (1,-0.1) {$1$};
    \node at (-0.1,1) {$1$};

    \coordinate (O) at (0,0);
    \coordinate (A) at (0, 1);
    \coordinate (B) at (0.7,1);
    \coordinate (Bl) at (0.55,1);
    \coordinate (B proj) at (0.7,0);
    \coordinate (C) at (1, 0.7); 
    \coordinate (C proj) at (0, 0.7); 
    \coordinate (D) at (1,0); 

    \coordinate (A') at (0, 0.99);
    \coordinate (B') at (0.7,0.99);
    \coordinate (C') at (0.99, 0.7); 
    \coordinate (D') at (0.99,0); 

\coordinate (An) at (0, 1);
    \coordinate (Bn) at (0.85,1);
    \coordinate (Bnl) at (1,0.99);
    \coordinate (Cn) at (1, 0.85); 
    \coordinate (Cnl) at (1, 0.85); 
    \coordinate (Dn) at (1,0); 

    \fill[gray!17] (O) -- (An) -- (Bn) -- (Cn) -- (Dn) -- (O) -- cycle;
    \fill[gray!60] (O) -- (A) -- (B) -- (C) -- (D) -- cycle;


\draw[thin, gray! 40] (0.7,0) -- (0.7,1);
\draw[thin, gray! 40] (0,0.7) -- (1,0.7);

    \draw[thick] (O) -- (0,0.99) node[midway, above left, black, font=\tiny] {}; 
    \draw[densely dotted, thick] (0.01, 1) -- (0.84,1) node[midway, right, black, font=\tiny] {};
    \draw[densely dotted, thick] (B) -- (C) node[midway, below right, black, font=\tiny] {};
     \draw[densely dotted, thick] (Bn) -- (Cn) node[midway, below right, black, font=\tiny] {};
     \draw[densely dotted, thick] (Cn) -- (C) node[midway, right, black, font=\tiny] {};
    \draw[densely dotted, thick] (1, 0.69) -- (1,0.01) node[midway, right, black, font=\tiny] {};
     \draw[thick] (0.99,0) -- (O) node[midway, right, black, font=\tiny] {};
    
     \node[above right] at (A) {$\Resize{0.2cm}{A}$};
     \node[above] at (Bl) {$ \Resize{1.5cm}{  B=\textstyle (\frac{2n-1}{2n}, 1) } $};
     \node[above] at (Bnl) {$\Resize{1.5cm}{ B'=\textstyle (\frac{2n}{2n+1}, 1)} $};
     \node[right] at (C) {$ \Resize{1.5cm}{ C=\textstyle (1,\frac{2n-1}{2n}) } $};
     \node[right] at (Cnl) {$ \Resize{1.5cm}{  C'=\textstyle (1,\frac{2n}{2n+1})}$};
      \node[above right] at (D) {$\Resize{.2cm}{D}$};
    \node[below left] at (O) {$\Resize{0.25cm}{O}$};
    
    \fill (O) circle (0.3pt);
    \draw (A) circle (0.3pt);
    \draw (B) circle (0.3pt);
    \draw (C) circle (0.3pt);
    \draw (Bn) circle (0.3pt);
    \draw (Cn) circle (0.3pt);
    \draw (D) circle (0.3pt);
    
    \node at (0.50, 0.50) {$\mathcal{P}$};
\end{tikzpicture}
\caption{$\mathcal{S}_{\mathrm{full}}$-boundedness region $\mathcal{P}$ which is the union of open pentagon with vertices $O,A,B,C,D$ and the line segments $[O,A)$ and $[O,D)$.}
\label{fig: full max region}
\end{figure}

\noindent \textit{The bilinear Kor\'anyi spherical maximal function.} Let $\mathfrak{G}$ denotes the surface 
\begin{equation}\label{unit sph in prod Hn}
    \mathfrak{G}:= \left\lbrace (z_{1}, t_{1}, z_{2}, t_{2}) \in \Ha \times \Ha : \Phi(z_{1}, t_{1}, z_{2}, t_{2})=0 \right\rbrace,
\end{equation}
where $\Phi$ is the smooth function $\R^{4n+2}$ defined as
\begin{equation*}
    \Phi(z_{1}, t_{1}, z_{2}, t_{2}) := (\|z_{1}\|^2 + \|z_{2} \|^{2})^2 + t_{1}^2 + t_{2}^2-1.
\end{equation*}
Let $f, g\in C^{\infty}_{c}(\Ha)$ and $r>0.$ We consider the bilinear Kor\'anyi averages as
\begin{equation}
\label{def:single:scale}
S_{r}(f,g)(x):= \int_{\mathfrak{G}} f \Big( x \cdot \delta_{r}(z_{1},t_{1}) \Big) \, g \Big( x \cdot \delta_{r}(z_{2},t_{2}) \Big) d \sigma(z_{1}, z_{2}, t_{1}, t_{2}),\quad x\in \Ha,
\end{equation}
where $d\sigma$ is the normalized induced Lebesgue measure on $\mathfrak{G}$ in $\R^{4n+2}.$ The averaging operator $S_{r}$ respects the scaling property below, which will be frequently used without mentioning
\begin{equation}
    S_{r}(f,g) = S_{1}(f \circ \delta_{r}, g \circ \delta_{r}) \circ \delta_{1/r},\quad 0<r<\infty.
    \label{dil}
\end{equation}
For brevity, we denote $S_1$ by $S.$ The bilinear Kor\'anyi spherical maximal function is then defined as 
\begin{align}
\label{def:bilinear:koranyi}
\mathcal{S}_{\mathrm{full}}(f, g)(x):=\sup_{r>0}|S_{r}(f, g)(x)|,\quad x\in \Ha,   
\end{align}
for $f, g\in C_{c}^{\infty}(\Ha).$ Now we state our main result. 
\begin{theorem}
\label{full:main}
Let $n\geq 2.$ Suppose $\big(\frac{1}{p_1}, \frac{1}{p_2}\big)$ is contained in the region $\mathcal{P}$, which consists of the open pentagon  with the corners $O=(0,0)$, $A=(0,1)$, $B=(\frac{2n-1}{2n},1)$, $C=(1,\frac{2n-1}{2n})$ and $D=(1,0)$, together with half-open line segments $[O,A)$ and $[O,D),$   see Figure \ref{fig: full max region}.  Let $0<p \leq \infty$ obeying the H\"older relation $\frac{1}{p} = \frac{1}{p_1} + \frac{1}{p_2}$. Then, we have 
\begin{equation}
   \| \mathcal{S}_{\mathrm{full}}(f,g)\|_{L^{p}(\Ha)} \leq C \|f\|_{L^{p_{1}}(\Ha)} \|g\|_{L^{p_{2}}(\Ha)}.
   \label{full:max}
\end{equation}
\end{theorem}

The reason for calling \eqref{def:bilinear:koranyi} as the bilinear analogue of the linear Kor\'anyi maximal operator is as follows. The product of two Heisenberg groups may be realized as $\R^{4n}\times \R^2,$ endowed with a step-two M\'etivier group structure with vertical part $\R^2.$ Now, the Carnot--Caratheodory distance $d((z_1, z_2, t_1, t_2))$ between the point $(z_1, z_2, t_1, t_2)$ and the identity $\vec{0}$ is homogeneous of degree one with respect to the scaling $\tilde{\delta}_{r}((z_1, z_2, t_1, t_2)))=(rz_1, r z_2, r^2 t_1, r^2\,t_2),$ further, it is well-known (see \cite{VaroBook}) that $d(z_1, z_2, t_1, t_2)\sim \big((\|z_{1}\|^2 + \|z_{2} \|^{2})^2 + t_{1}^2 + t_{2}^2\big)^{1/4},$ and this dictates the surface considered in \eqref{def:single:scale}.

Theorem \ref{full:main} gives a sufficient range of exponents for the boundedness of $\mathcal S_{\mathrm{full}}$. Our next result provides a necessary condition, showing that there remains a gap between the sufficient and necessary ranges, see Figure \ref{fig: full max region}. 

\begin{prop}
\label{nece:full:new}
Let $n\geq 1.$ Let $1\leq p_1, p_2\leq \infty$ and $0<p<\infty.$ If we have \[\|\mathcal{S}_{\rm{full}}(f, g)\|_{L^p(\Ha)}\leq C\|f\|_{L^{p_1}(\Ha)}\|g\|_{L^{p_2}(\Ha)},\] then we must have
\begin{align}
\label{Nece:2}
&\frac{1}{p} = \frac{1}{p_1}+\frac{1}{p_2}\leq 2-\frac{1}{2n+1}.
\end{align}
That is, $(\frac{1}{p_1}, \frac{1}{p_2})$ must be contained in the closed pentagonal region with corners $O=(0,0)$, $A=(0,1)$, $B'=(\frac{2n}{2n+1},1)$, $C'=(1,\frac{2n}{2n+1})$ and $D=(1,0)$, see Figure \ref{fig: full max region}. 
\end{prop}
One observes that the H\"older relation $\frac{1}{p}=\frac{1}{p_1}+\frac{1}{p_2}$ in Proposition \ref{nece:full:new} is a consequence of the scaling law \eqref{dil}.

\noindent\textit{Previous Work.} To put our results in the context, recall that in recent years the study of Euclidean bilinear spherical maximal function, defined as
$$\mathcal{M}(f, g)(x):=\sup_{r>0}\bigg|\int_{S^{2n-1}}f(x-ry)g(x-rz)\, d\sigma_{2n-1}(y, z)\bigg|,$$ 
has attracted a great deal of research, starting with the work of Geba \textit{et al} in \cite{IoseMRL}, Barrionuevo \textit{et al} \cite{GraMRL} and subsequent improvements in \cite{GHH21} and \cite{HHY20}. Finally, in \cite{LeeJFA}, Jeong and Lee obtained a sharp range of exponents for the Lebesgue space boundedness of $\mathcal{M}$ when $n\geq 2$. They proved that $\mathcal{M}$ is bounded from $L^p(\R^n)\times L^q(\R^n)$ to $L^r(\R^n)$ if and only if $r>\frac{n}{2n-1},$ $\frac{1}{r}=\frac{1}{p}+\frac{1}{q}, 1\leq p, q\leq \infty,$ except the exponent $(p,q,r)= (1, \infty, 1)$ or $ (\infty, 1, 1)$, where they obtained corresponding weak-type estimates. The case for $n=1$ is independently addressed by Christ--Zhou \cite{CZ24}, and Dosidis--Ramos \cite{DR24} and their result states that $\mathcal{M}: L^p(\R)\times L^{q}(\R)\to L^{r}(\R)$ if and only if $\frac{1}{r}=\frac{1}{p}+\frac{1}{q},$ $2<p, q\leq \infty.$ For the lacunary maximal function
$$\mathcal{M}_{\rm{lac}}(f, g)(x)=:=\sup_{k\in \mathbb{Z}}\bigg|\int_{S^{2n-1}}f(x-2^{-k}y)g(x-2^{-k}z)\, d\sigma_{2n-1}(y, z)\bigg|,$$
we refer to \cite{BF24} for $n\geq 2$ and to \cite{CZ24} for the case $n=1.$ For sparse bounds of the bilinear spherical maximal function, we refer to \cite{BFO23}, while \cite{IoseJGA} contains Lebesgue space estimates for the single-scale operators.

On the Heisenberg group, as already mentioned, Cowling \cite{Cowling} extended Stein's theorem for the Kor\'anyi spherical maximal as stated in Theorem~\ref{Koranyi:Thm}. In recent times, in \cite{PT}, Ganguly--Thangavelu studied the lacunary maximal function associated with Kor\'anyi spherical averages, $A_{r}^{K}$, and subsequently in \cite{Rajula} the author obtained $L^p$ improving estimates for the local Kor\'anyi maximal function and obtained sparse operator bounds and weighted estimates. Besides the Kor\'anyi spherical maximal function, there is another important maximal operator on $\Ha$, defined by averaging over the Euclidean sphere $S^{2n-1}\times \{0\}\subset \C^n\times \R.$ More precisely,
$$M_{S}f(x)=\sup_{r>0}|A_{r}f(x)|,$$
where $A_r$ stands for Nevo--Thangavelu spherical means given by 
\begin{align}
\label{NT}
A_{r} f(x):=\int_{S^{2n-1}} { \textstyle f\big(\underline{x}-r \underline{u}, \overline{x}- {\textstyle\frac{1}{2} } r  \Im(\underline{x} \cdot \bar{\underline{u}})\big)}\,d \sigma_{2n-1}(\underline{u}),\quad x=(\underline{x}, \overline{x}) \in \mathbb{H}^n.    
\end{align}
The operator $M_{S}$ will also play a pivotal role in our proofs. It was first studied by Nevo--Thangavelu in \cite{NeT}. It is well-known, in  the case of $n\geq 2$, that $M_{S}$ is bounded on $L^p(\mathbb{H}^n)$ if and only if $p > \frac{2n}{2n-1}$ which was independently proved by M\"uller--Seeger in \cite{MS} and by Narayanan--Thangavelu in \cite{NaT}. Very recently, this result is extended to all two-step nilpotent Lie groups by Ryu--Segger in \cite{SeegerAdv}. Moreover, in \cite{BHRT, Roos} $L^p$ improving estimates for these averages and sparse operator bounds are established. We also refer the recent work \cite{HickmanJFA} where the authors developed a general framework to study lacunary maximal functions on homogeneous groups.

When it comes to the bilinear setting on $\Ha$, very recently, in \cite{Ghosh-Rajesh}, we studied the bilinear Nevo--Thangavelu maximal function and established its $L^{p_1}(\Ha)\times L^{p_2}(\Ha)\to L^p(\Ha)$ boundedness. The present article is devoted to obtaining the corresponding theory for the bilinear Kor\'anyi spherical maximal function. Unlike the Euclidean setting, the presence of the symplectic form in the group law changes geometry of the surface with the centre, requiring new ideas to analyze the associated bilinear averaging operators.\\

\noindent\textit{Ideas behind the proofs.} We conclude this section with a few remarks on the proof of Theorem~\ref{full:main}. Our argument is inspired by the slicing method introduced by Jeong and Lee \cite{LeeJFA} in the Euclidean setting. The main idea is to decompose the bilinear Kor\'anyi averages into lower-dimensional averaging operators, which naturally leads to several auxiliary maximal operators. We now recall the ones that play a crucial role in our proof. The first one is the ergodic maximal operator $\Lambda$, defined by
\begin{equation}
\label{ergodic-max}
\Lambda g(x)
:=
\sup_{r>0}
\bigg|
\frac{1}{r}
\int_{0}^{r}
A_{s}g(x)\,ds
\bigg|,
\qquad x\in\Ha.
\end{equation}
By Birkhoff's maximal theorem, it is well known that $\Lambda$ is bounded on $L^q(\Ha)$ for every $1<q\leq\infty$; see Proposition~3.6.1 of \cite{ThangBook}.

The second auxiliary operator is a one-dimensional Euclidean bilinear maximal operator. For functions $h_1,h_2$ on $\R$, we define
\begin{equation}
\label{max:avg:bi}
\mathcal{B}(h_1,h_2)(\overline{x})
:=
\sup_{r>0}
\bigg|
\frac{1}{r}
\int_{0}^{r}
B_{\rho}(h_1,h_2)(\overline{x})\,d\rho
\bigg|,
\qquad
\overline{x}\in\R,
\end{equation}
where
\[
B_{\rho}(h_1,h_2)(\overline{x})
=
\int_{0}^{2\pi}
h_1(\overline{x}+\rho\cos\theta)\,
h_2(\overline{x}+\rho\sin\theta)\,
d\theta.
\]
The operator $\mathcal{B}$ arises naturally from the decomposition of the bilinear Kor\'anyi averages and plays a central role in our analysis, more precisely, we show that 
\begin{align}
\label{main:pointwise}
 \mathcal{S}_{\rm{full}}(f, g)(\underline{x}, \overline{x})\lesssim \mathcal{B}(\scriptsize M_{S}  f  \big( \underline{x}, \cdot \big),   \scriptsize \Lambda  g   \big( \underline{x}, \cdot \big))(\overline{x})=  \sup_{r>0} \frac{1}{r} \int_{0}^{r} B_{ \rho}  \Big( { \scriptsize M_{S}  f  \big( \underline{x}, \cdot \big)  }  \, , \,   { \scriptsize \  \Lambda  g   \big( \underline{x}, \cdot \big) } \Big)  (\overline{x})d \rho
\end{align}
for any $x=(\underline{x}, \overline{x})\in \Ha,$ and a similar estimate interchanging the role of $f$ and $g$ on the right hand side of \eqref{main:pointwise}. Here, of course, $B_{\rho}\, (\text{and}\,\, \mathcal{B})$ acts on the central direction $\overline{x}.$

Finally, it remains to establish the boundedness of the auxiliary operator $\mathcal{B}$. This is the content of the following theorem. 
\begin{theorem}
\label{euclidean:avg:thm}
Let $1<p_1, p_2\leq \infty$ and $\frac{1}{p}=\frac{1}{p_1}+\frac{1}{p_2}.$ Then $\mathcal{B}$ maps $L^{p_1}(\mathbb{R})\times L^{p_2}(\mathbb{R})$ boundedly to $L^{p}(\mathbb{R}).$ If either $p_1$ or $p_2$ is $1,$ corresponding weak-type estimates hold.   
\end{theorem}

\begin{remark}
The operator $\mathcal{B}$ plays a crucial role in our proof. Indeed, one could instead use the bilinear spherical maximal function
\[
\mathcal{M}(f,g)(\overline{x})
:=
\sup_{\rho>0}
\bigg|
\int_{0}^{2\pi}
f(\overline{x}+\rho\cos\theta)\,
g(\overline{x}+\rho\sin\theta)\,
d\theta
\bigg|,
\qquad
\overline{x}\in\R,
\]
in the pointwise estimate \eqref{main:pointwise}. However, it is known that $\mathcal{M}$ is bounded from $L^{p_1}(\R)\times L^{p_2}(\R)$ to $L^p(\R)$ if and only if
\[
\frac{1}{p}
=
\frac{1}{p_1}
+
\frac{1}{p_2},
\qquad
2<p_1,p_2\leq\infty.
\]
Consequently, replacing $\mathcal{M}$ by the averaging operator $\mathcal{B}$ in \eqref{main:pointwise} is one of the key ingredients that allows us to obtain a substantially larger range of exponents in Theorem~\ref{full:main}.
\end{remark}

The article is organized as follows. In Section~\ref{sec:slice} we prove  Theorem~\ref{euclidean:avg:thm} and Theorem~\ref{full:main}. In section~\ref{sec:nece} we prove  Proposition~\ref{nece:full:new}. Finally, in Section~\ref{multi}, we obtain multilinear generalizations of our results. The notation $A\lesssim B$ represents that $A\leq C B,$ for some constant $C$ independent of $A, B$, and $A\sim B$ indicates $A \lesssim B$ and $B \lesssim A.$ The normalized surface measure on the unit sphere $S^{k-1}\subset \R^k$ is denoted by $d\sigma_{k-1}.$

\section{Slicing and the proof of Theorem~\ref{full:main}}
\label{sec:slice}

In this section, we present the proof our main result Theorem~\ref{full:main}. The proof is divided into several parts. 

\subsection{An Euclidean maximal operator}
We start by proving Theorem~\ref{euclidean:avg:thm}. First, let us recall one of the fundamental operators in the study of the bilinear (or multilinear) Calder\'on-Zygmund theory, namely the bilinear Hardy-Littlewood maximal operator on $\R$ defined as
\[\mathcal{M}_{\rm{HL}}(h_1, h_2)(\overline{x})=\sup_{I\ni \overline{x}} \frac{1}{|I|}\int_{I}|h_1(y)|\, dy \frac{1}{|I|}\int_{I}|h_2(z)|\, dz,\quad \overline{x}\in \R,\]
where the supremum is taken over all intervals $I$ containing $\overline{x}.$ Using H\"older's inequality and $\mathcal{M}_{\rm{HL}}(h_1, h_2)(\overline{x})\leq M_{\rm{HL}}(h_1)(\overline{x})\,M_{\rm{HL}}(h_{2})(\overline{x}),$ where $M_{\rm{HL}}$ is the standard uncentered (linear) Hardy-Littlewood maximal operator, it is easy to see that $\mathcal{M}_{\rm{HL}}: L^{p_1}(\R)\times L^{p_2}(\R)\to L^{p}(\R)$ for all $1<p_{i}<\infty$ and $\frac{1}{p}=\frac{1}{p_1}+\frac{1}{p_2}.$ In fact, $\mathcal{M}_{\rm{HL}}$ is strictly smaller than the product $M_{\rm{HL}}h_{1} M_{\rm{HL}}h_{2},$ and played a pervasive role in the theory of multilinear $A_{\vec{P}}$ weights, see \cite{Lerner-adv}.

\begin{proof}[Proof of Theorem~\ref{euclidean:avg:thm}]

 Fix $r>0.$ Without loss of generality we assume $h_1, h_2\geq0.$ Observe that
\begin{align}
\nonumber&\frac{1}{r}\int_{0}^{r}B_{\rho}(h_1, h_2)(\overline{x})\, d\rho\\
\nonumber& =  \frac{1}{r}\int_{0}^{r} \int_{0}^{2\pi} h_{1}(\overline{x}+\rho \cos \theta) h_{2}(\overline{x}+\rho \sin \theta)\, d\theta\, d\rho\\
\nonumber &=\frac{1}{r}\int_{0}^{r} \int_{0}^{2\pi} h_{1}(\overline{x}+\rho \cos \theta) h_{2}(\overline{x}+\rho \sin \theta)\, \frac{1}{\rho} \rho\, d\rho\, d\theta.
\end{align}
Using cartesian coordinates, the above expression turns into
\[\frac{1}{r}\int\int_{\sqrt{u^2+v^2}\leq r} h_{1}(\overline{x}+u) h_{2}(\overline{x}+v)\, \frac{1}{\sqrt{u^2+v^2}} du\,dv.\]
For $k\in \N\cup \{0\},$ we denote $R_k = \{(u, v) \in \mathbb{R}^2 : 2^{-k-1}r < \sqrt{u^2 + v^2} \leq 2^{-k}r \},$ and decompose the above integral as 
$$\frac{1}{r}\sum_{k\geq 0} \int\int_{R_k} h_{1}(\overline{x}+u) h_{2}(\overline{x}+v)\, \frac{1}{\sqrt{u^2+v^2}} du\,dv.$$
Subsequently, 
\begin{align}
\nonumber&\big|\frac{1}{r}\int_{0}^{r}B_{\rho}(h_1, h_2)(\overline{x})\, d\rho\big|\\
\nonumber &\lesssim \frac{1}{r}\sum_{k\geq 0} \frac{1}{2^{-k-1}r} \int\int_{R_k} h_{1}(\overline{x}+u)\, h_{2}(\overline{x}+v)|\, du\,dv\\
\nonumber &\lesssim \frac{1}{r}\sum_{k\geq 0} 2^{-k}r \frac{1}{(2^{-k}r)^2} \int_{|u|\leq 2^{-k}r} \int_{|v|\leq 2^{-k}r} h_{1}(\overline{x}+u) \,h_{2}(\overline{x}+v)\, du\,dv\\
&\lesssim \mathcal{M}_{\rm{HL}}(h_1, h_2)(\overline{x})\sum_{k\geq 0} 2^{-k}\lesssim \mathcal{M}_{\rm{HL}}(h_1, h_2)(\overline{x}).
\end{align}
Taking supremum over $r>0$ we obtain that 
\begin{align}
\label{max:avg:bi}
\mathcal{B}(h_1, h_2)(\overline{x})=\sup_{r>0}\big|\frac{1}{r}\int_{0}^{r}B_{\rho}(h_1, h_2)(\overline{x})\, d\rho\big|\lesssim \mathcal{M}_{\rm{HL}}(h_1, h_2)(\overline{x}).
\end{align}
Now the result follows since $\mathcal{M}_{\rm{HL}}: L^{p_1}(\R)\times L^{p_2}(\R)\to L^{p}(\R)$ boundedly for all $1<p_{i}<\infty$ and $\frac{1}{p}=\frac{1}{p_1}+\frac{1}{p_2}.$ If $1\leq p_1, p_2\leq \infty$ with either $p_1=1$ or $p_2=1,$ $\mathcal{B}$ maps $L^{p_1}(\R)\times L^{p_2}(\R)$ to $L^{p, \infty}(\R)$ since the same holds for $\mathcal{M}_{\rm{HL}}.$  
\end{proof}

\begin{remark}
We make a discursive remark regarding the weighted estimates for $\mathcal{B}.$ Let $1<p_1, p_2<\infty,$ set $\vec{P}=(p_1, p_2).$ Denote $v_{\vec{w}}=w_{1}^{p/p_1}w_{2}^{p/p_2}.$ A tuple of weights $\vec{w}=(w_1, w_2)\in A_{\vec{P}}$ if and only if
$$\sup_{I} \left( \frac{1}{\vert{}I\vert{}} \int_I v_{\vec{w}} \, dx \right)^{\frac{1}{p}} \prod_{i=1}^2 \left( \frac{1}{\vert{}I\vert{}} \int_I w_i(x)^{1-p_i'} \, dx \right)^{\frac{1}{p_i'}} < \infty,$$
where the supremum ranges over all intervals in $\mathbb{R}.$ Then, as an application of our Theorem~\ref{euclidean:avg:thm} and Theorem 3.7 in \cite{Lerner-adv}, for $1<p_1, p_2<\infty,$ $\frac{1}{p}=\frac{1}{p_1}+\frac{1}{p_2},$ and $\vec{w}\in A_{\vec{P}}$ the weighted estimate $$\|\mathcal{B}(h_1, h_2)\|_{L^p(v_{\vec{w}})}\leq C \|h_1\|_{L^{p_1}(w_1)} \|h_2\|_{L^{p_2}(w_2)}$$
holds true.
\end{remark}

\subsection{Slicing and domination}
Our next proposition is the key result of this article. This uses the slicing technique as in the paper \cite{LeeJFA}, however, since we are in the Heisenberg group, we need to deal with additional difficulties arising from the group structure. For any continuous function $G$ on $\R^{2n}, n\geq 2,$ we recall the slicing formula from \cite{LeeJFA}
\begin{align}
\label{slicing-formula-Lee}
&\int_{S^{2n-1}} G(z_1, z_2)\,d\sigma_{2n-1}(z_1, z_2)\\
\nonumber&= \int_{B^{n}(0, 1)}\int_{S^{n-1}} G(z_1, \sqrt{1-\|z_1\|^2} z_2)\, d\sigma_{n-1}(z_2)(1-\|z_1\|^2)^{(n-2)/2}\, dz_1,   
\end{align}
as this will be used directly in our proof.

\begin{prop} Let $f, g\in C_{c}^{\infty}(\Ha).$ The following pointwise estimate holds
\begin{align}
 \mathcal{S}_{\rm{full}}(f, g)(\underline{x}, \overline{x})\lesssim \mathcal{B}(\scriptsize M_{S}  f  \big( \underline{x}, \cdot \big),   \scriptsize \Lambda  g   \big( \underline{x}, \cdot \big))(\overline{x})=  \sup_{r>0} \frac{1}{r} \int_{0}^{r} B_{ \rho}  \Big( { \scriptsize M_{S}  f  \big( \underline{x}, \cdot \big)  }  \, , \,   { \scriptsize \  \Lambda  g   \big( \underline{x}, \cdot \big) } \Big)(\overline{x})\, d\rho,
 \label{slice:eqn}
\end{align}
for any $x=(\underline{x}, \overline{x})\in \Ha.$ Also, a similar estimate holds interchanging the role of $f$ and $g$ on the right hand side of \eqref{slice:eqn}.       
\end{prop}

\begin{proof}
Motivated by \cite{LeeJFA} we use the slicing argument to express the averaging operators $S_{r}$ in terms of linear horizontal averaging operators $A_{s}$. To start with, we recall the following formula. Let $\Omega\subset \mathbb{R}^{N}$ be a $(N-1)$ dimensional hypersurface given by $\Omega:=\{\theta\in \mathbb{R}^N: \Psi(\theta)=0\},$ with $\|\nabla \Psi(\theta) \|\neq 0$ on $\Omega.$ Then
\begin{align}
\label{slicing:formula}
\int_{\Omega} G(\theta)\frac{d\mu(\theta)}{\|\nabla \Psi(\theta)\|}=\int_{\mathbb{R}^{N}} G(x)\delta(\Psi)\, dx,
\end{align}
where $\mu$ is the induced Lebesgue measure on $\Omega$ and $G\in C(\mathbb{R}^N).$ Without loss of generality we assume $f, g\geq 0.$ Recall $$\mathfrak{G}:= \left\lbrace (z_{1}, t_{1}, z_{2}, t_{2}) \in \R^{2n+1} \times \R^{2n+1} : \Phi(z_{1}, t_{1}, z_{2}, t_{2})=0 \right\rbrace,$$ where
\begin{equation*}
    \Phi(z_{1}, t_{1}, z_{2}, t_{2}) := (\|z_{1}\|^2 + \|z_{2} \|^{2})^2 + t_{1}^2 + t_{2}^2 - 1.
\end{equation*}
To apply the formula \eqref{slicing:formula}, we compute 
$$\nabla_{(z_{1}, z_{2})} \Phi (z_{1}, z_{2},t_{1},t_{2}) = 4 (\| z_{1} \|^2 + \|z_{2} \|^2) (z_{1}, z_{2})$$
and
$$\nabla_{(t_{1}, t_{2})} \Phi (z_{1}, z_{2},t_{1},t_{2}) = 2 ( t_{1}  , t_{2}).$$
For simplicity, we denote $\vec{\mathfrak{z}}:=(z_1, z_2)\in \mathbb{C}^{2n}$ and $\vec{\mathfrak{t}}=(t_1, t_2),$ and thus \begin{equation*}
\begin{split}
        \| \nabla \Phi \| & = \Big\| \Big( 4 (\| z_{1} \|^2 + \|z_{2} \|^2) (z_{1}, z_{2}) \, , \, 2 ( t_{1}  , t_{2})  \Big) \Big\| \\
        & =  \Big( 16 \| \vec{\mathfrak{z}} \|^6 + 4 \| \vec{\mathfrak{t}} \|^2 \Big)^{1/2} \\
        & = 2 \Big( 4( 1-  \| \vec{\mathfrak{t}} \|^2 )^{3/2} +  \| \vec{\mathfrak{t}} \|^2 \Big)^{1/2}\, \quad\text{on}\quad \mathfrak{G}.
\end{split}
\end{equation*}
It is easy to see that there exist two constants $c_1, c_2>0$ such that $c_{1}\leq \| \nabla \Phi \|=2 \Big( 4( 1-  \| \vec{\mathfrak{t}} \|^2 )^{3/2} +  \| \vec{\mathfrak{t}} \|^2 \Big)^{1/2}\leq c_{2}$ holds on $\mathfrak{G}.$ Now applying the formula \eqref{slicing:formula} we derive that
\begin{equation*}
\begin{split}
&\int_{\mathfrak{G}} f \Big( x \cdot (z_{1},t_{1}) \Big) \, g \Big( x \cdot (z_{2},t_{2}) \Big) \frac{d \sigma(z_{1}, z_{2}, t_{1}, t_{2})}{\|\nabla \Phi(z_1, t_1, z_2, t_2)\|}\\
&= \int_{ (\vec{\mathfrak{z}}, \vec{\mathfrak{t}}) \in \R^{4n} \times \R^2}  \textstyle   f \Big( x \cdot (z_{1},t_{1}) \Big) \, g \Big( x \cdot (z_{2},t_{2}) \Big)  \ \delta (\Phi) \  d \vec{\mathfrak{z}} \, d \vec{\mathfrak{t}}.
\end{split}    
\end{equation*}
As $\|\nabla \Phi\|\sim 1,$ we obtain that
$$S(f, g)(x):=S_{1}(f, g)(x)\sim \int_{ (\vec{\mathfrak{z}}, \vec{\mathfrak{t}}) \in \R^{4n} \times \R^2}  \textstyle   f \Big( x \cdot (z_{1},t_{1}) \Big) \, g \Big( x \cdot (z_{2},t_{2}) \Big)  \ \delta (\Phi) \  d \vec{\mathfrak{z}} \, d \vec{\mathfrak{t}}.$$
For a fixed 
$\vec{\mathfrak{t}}=(t_1, t_2)$ with $\|\vec{\mathfrak{t}}\|< 1,$ a further application of the slicing argument \eqref{slicing:formula}  yields
$$S(f, g)(x)\sim \int_{\| \vec{\mathfrak{t}} \| <1} \int_{ \vec{\mathfrak{z}} \in \Omega_{\vec{\mathfrak{t}}} }  \textstyle   f \Big( x \cdot (z_{1},t_{1}) \Big) \, g \Big( x \cdot (z_{2},t_{2}) \Big) \frac{ d \nu_{\vec{\mathfrak{t}}} (\vec{\mathfrak{z}})}{ \| \big( \nabla_{\vec{\mathfrak{z}}} \Phi_{\vec{\mathfrak{t}}} \big) (\vec{\mathfrak{z}}) \|   }  \, d \vec{\mathfrak{t}}\,,$$
where $$\Omega_{\vec{\mathfrak{t}}}:=\{\vec{\mathfrak{z}}=(z_1, z_2)\in \mathbb{R}^{2n}: \Phi_{\vec{\mathfrak{t}}}(z_1, z_2):= (\|z_{1}\|^2+\|z_{2}\|^2)^2-(1-\|\vec{\mathfrak{t}}\|^2)= 0\},$$ and $d \nu_{\vec{\mathfrak{t}}}$ is the induced measure on $\Omega_{\vec{\mathfrak{t}}}$ in $\R^{2n}.$ Since $\Omega_{\vec{\mathfrak{t}}}$ is the $(4n-1)$-dimensional Euclidean sphere of radius $R:= (1 - \| \vec{\mathfrak{t}} \|^2)^{1/4},$ a further change of variables simplifies the above to
\begin{equation*}
\begin{split}
&S(f,g)(x)\\
&\sim \int_{ \| \vec{\mathfrak{t}} \| <1} \int_{ \vec{\mathfrak{z}} \in S^{4n-1} }  \textstyle   f \Big( x \cdot ( (1- \| \vec{\mathfrak{t}}\|^2)^{1/4}z_{1},t_{1}) \Big) \, g \Big( x \cdot ( (1- \| \vec{\mathfrak{t}}\|^2)^{1/4}z_{2},t_{2}) \Big) \frac{R^{4n-1}   d \sigma_{4n-1}(\vec{\mathfrak{z}})}{4R^3}\, d \vec{\mathfrak{t}}\\
&\sim \int_{ \|\vec{\mathfrak{t}} \| <1}(1 - \| \vec{\mathfrak{t}} \|^2)^{n-1} \int_{ \vec{\mathfrak{z} }\in S^{4n-1} }  \textstyle   f \Big( x \cdot ( (1- \| \vec{\mathfrak{t}}\|^2)^{1/4}z_{1},t_{1}) \Big) \, g \Big( x \cdot ( (1- \| \vec{\mathfrak{t}}\|^2)^{1/4}z_{2},t_{2}) \Big)d \sigma_{4n-1}(\vec{\mathfrak{z}})\, d \vec{\mathfrak{t}}.
\end{split}    
\end{equation*}
Using polar decomposition in the $\vec{\mf{t}}=(\rho \cos \theta, \rho \sin \theta)$ variable the above expression becomes

$$S(f,g)(x) \sim \int_{\rho=0}^{1}  {\scriptsize \rho \,  (1 - \rho^2)^{n-1}} I_{f, g}(\rho)\, d\rho,$$
where 
\begin{align*}
I_{f, g}(\rho):= \int_{ \theta=0}^{2 \pi}   \int_{\vec{\mathfrak{z}} \in S^{4n-1} }  {\scriptsize  f ( x . ( (1- \rho^2)^{1/4}z_{1}, \rho \cos \theta) ) } \, {\scriptsize g ( x . ( (1- \rho^2)^{1/4}z_{2}, \rho \sin \theta  ) ) }   d \sigma_{4n-1}(\vec{\mathfrak{z}})\,d \theta. 
\end{align*}

A further application of the slicing formula \eqref{slicing-formula-Lee} yields
\begin{equation}\label{Ifg}
 \begin{split}
&I_{f, g}(\rho)\\
&= \int_{ \theta=0}^{2 \pi}   \int_{ \mathfrak{z} \in S^{4n-1} }  {\scriptsize  f ( x . ( (1- \rho^2)^{1/4}z_{1}, \rho \cos \theta) ) } \, {\scriptsize g ( x . ( (1- \rho^2)^{1/4}z_{2}, \rho \sin \theta  ) ) }   d \sigma_{4n-1}(\vec{\mathfrak{z}}) d\theta  \\
& = \int_{ \theta=0}^{2 \pi}   \int_{z_{1} \in B^{2n}(0,1)}   {\scriptsize  f ( x . ( (1- \rho^2)^{1/4}z_{1}, \rho \cos \theta) ) } \\
& \ \ \ \ \ \ \ \ \ \ \ \ \ \ \ \int_{z_{2} \in S^{2n-1}} {\scriptsize g ( x . ( (1- \rho^2)^{1/4} {\scriptsize \sqrt{1-\|z_{1}\|^2}} z_{2}, \rho \sin \theta  ) ) }   d \sigma_{2n-1}(z_{2}) (1-\|z_{1}\|^2)^{n-1} d z_{1}d\theta.  
 \end{split}
\end{equation}
Now, the inner integral can be expressed as
\begin{equation*}
    \begin{split}
& \int_{z_{2} \in S^{2n-1}} {\scriptsize g ( x . ( (1- \rho^2)^{1/4} {\scriptsize \sqrt{1-\|z_{1}\|^2}} z_{2}, \rho \sin \theta  ) ) }   d \sigma_{2n-1}(z_{2}) \\
& = \int_{z_{2} \in S^{2n-1}} {\scriptsize g \left( x \cdot (\mathbf{0}, \rho \sin \theta) \cdot \delta_{ (1- \rho^2)^{1/4}  \sqrt{1-\|z_{1}\|^2}  } (  z_{2}, 0 )  \right) }  d \sigma_{2n-1}(z_{2}) \\
& = \big( A_{ (1- \rho^2)^{1/4}  \sqrt{1-\|z_{1}\|^2} } \, g \big) \big( x \cdot (\mathbf{0}, \rho \sin \theta) \big).
    \end{split}
\end{equation*}
Therefore, \eqref{Ifg} in polar coordinates becomes
\begin{equation*}
\begin{split} 
&I_{f, g}(\rho)\\
& =\int_{ \theta=0}^{2 \pi}   \int_{z_{1} \in B^{2n}(0,1)}   {\scriptsize  f ( x . ( (1- \rho^2)^{1/4}z_{1}, \rho \cos \theta) ) } \ { \scriptsize \big( A_{ (1- \rho^2)^{1/4}  \sqrt{1-\|z_{1}\|^2} } \, g \big) \big( x \cdot (\mathbf{0}, \rho \sin \theta) \big)    (1-\|z_{1}\|^2)^{n-1} } d z_{1} d\theta\\
& = \int_{ \theta=0}^{2 \pi}   \int_{s=0}^{1} \int_{ \omega \in S^{2n-1}(0,1)}   {\scriptsize  f ( x . ( (1- \rho^2)^{1/4} s \omega, \rho \cos \theta) ) }\, d \sigma_{2n-1}(\omega) \,\\
& \ \ \ \ \ \ \ \ \ \ \ \ \ \ \ \ \  { \scriptsize \big( A_{ (1- \rho^2)^{1/4}  \sqrt{1- s^2} } \, g \big) \big( x \cdot (\mathbf{0}, \rho \sin \theta) \big)    \, (1-s^2)^{n-1} s^{2n-1} } ds\,d\theta \\
& = \int_{ \theta=0}^{2 \pi}   \int_{s=0}^{1} { \scriptsize \big( A_{ (1- \rho^2)^{1/4}  s } \, f \big) \big( x \cdot (\mathbf{0}, \rho \cos \theta) \big)  }   { \scriptsize \big( A_{ (1- \rho^2)^{1/4}  \sqrt{1- s^2} } \, g \big) \big( x \cdot (\mathbf{0}, \rho \sin \theta) \big) } (1-s^2)^{n-1} s^{2n-1} ds\,d\theta.\\
\end{split}    
\end{equation*}

Combining the above, if we write $x$ as $x= (\underline{x}, \overline{x}) \in \Ha $, then 
\begin{align}
\nonumber& S(f,g)(\underline{x}, \overline{x})\\
\nonumber& \sim  \int_{\rho=0}^{1} {\scriptsize \rho \,  (  1 - \rho^2    )^{n-1}  } \int_{ \theta=0}^{2 \pi}   \int_{z_{1} \in B^{2n}(0,1)}   {\scriptsize  f ( x . ( (1- \rho^2)^{1/4}z_{1}, \rho \cos \theta) ) } \\
& \ \ \ \ \ \ \ \ \ \ \ \ \ \ \ \int_{z_{2} \in S^{2n-1}} {\scriptsize g ( x . ( (1- \rho^2)^{1/4} {\scriptsize \sqrt{1-\|z_{1}\|^2}} z_{2}, \rho \sin \theta  ) ) }   d \sigma_{2n-1}(z_{2}) (1-\|z_{1}\|^2)^{n-1} d z_{1}\, d\theta\label{slice:final:1}\\ 
\nonumber&=\int_{\rho=0}^{1}  \int_{ \theta=0}^{2 \pi} \int_{s=0}^{1} {\scriptsize \rho \,  (  1 - \rho^2    )^{n-1}  }           { \scriptsize \Big( A_{ (1- \rho^2)^{1/4}  s } \, f \Big) \big( \underline{x}, \overline{x} +  \rho \cos \theta \big)  } \\
& \hspace{10.5em} { \scriptsize \Big( A_{ (1- \rho^2)^{1/4}  \sqrt{1- s^2} } \, g \Big) \big( \underline{x}, \overline{x} +  \rho \sin \theta \big) } \,  (1-s^2)^{n-1} s^{2n-1} ds       \,             d \rho \, d \theta \label{slice:final:2}.
\end{align}

The formulas \eqref{slice:final:1} and \eqref{slice:final:2} will be used frequently in most of our proofs. In the above expression \eqref{slice:final:2}, we dominate the horizontal average on the $f$-component by the Nevo--Thangavelu maximal function $M_{S}$ to arrive at the following
\begin{align}
\label{single scale dom:1}
\nonumber S(f,g) (\underline{x}, \overline{x}) & \lesssim \int_{\rho=0}^{1}  {\scriptsize \rho \,  (  1 - \rho^2    )^{n-1}  } \int_{ \theta=0}^{2 \pi} { \scriptsize M_{S}  f  \big( \underline{x}, \overline{x} +  \rho \cos \theta \big)  }\\
&   \hspace{7em}        \int_{s=0}^{1} { \scriptsize \Big( A_{ (1- \rho^2)^{1/4}  \sqrt{1- s^2} } \, g \Big) \big( \underline{x}, \overline{x} +  \rho \sin \theta \big) }   \,  (1-s^2)^{n-1} s^{2n-1} ds       \,             d \rho \, d \theta.   
\end{align}

Using first the change of variables $s \rightarrow \sqrt{ 1- (1- \rho^2)^{-1/2} s^2 }$ in \eqref{single scale dom:1} and then dominating polynomial factors in $s$ by some constants, we derive 

\begin{align}
\label{single scale dom}          
S(f,g) (\underline{x}, \overline{x}) & \lesssim \int_{\rho=0}^{1}  \rho \,  (  1 - \rho^2    )^{n-1}   \int_{\theta=0}^{2 \pi} { M_{S}  f  \big( \underline{x}, \overline{x} +  \rho \cos \theta \big)} \\         
\nonumber & \quad \frac{1}{(1- \rho^2)^{1/4}}    \int_{s=0}^{ (1- \rho^2)^{1/4} } ( A_{ s } \, g ) \big( \underline{x}, \overline{x} +  \rho \sin \theta \big) \,ds\, d \rho \, d \theta.
\end{align}
At this point, we invoke the operator $\Lambda$ in the $g$-component and obtain the following
\begin{align} 
\nonumber S(f,g) (\underline{x}, \overline{x})&\lesssim \int_{\rho=0}^{1}  {\scriptsize \rho \,  (  1 - \rho^2    )^{n-1}  }         \int_{ \theta=0}^{2 \pi} { \scriptsize M_{S}  f  \big( \underline{x}, \overline{x} +  \rho \cos \theta \big)  }    { \scriptsize \  \Lambda  g  \big( \underline{x}, \overline{x} +  \rho \sin \theta \big) } \, d \theta         \,             d \rho  \\
 & = \int_{\rho=0}^{1}  {\scriptsize \rho \,  (  1 - \rho^2    )^{n-1}  }         B_{\rho} \Big( { \scriptsize M_{S}  f  \big( \underline{x}, \cdot \big)  }  \, , \,   { \scriptsize \  \Lambda  g  \big( \underline{x}, \cdot \big) } \Big)   (\overline{x})\,d \rho, \label{final:single} 
\end{align}
where $B_{\rho}$ is the Euclidean bilinear spherical means acting on the central variable $\overline{x}$ given by
\begin{equation*}
 B_{\rho} \big(h_1, h_2 \big)(\overline{x})  := \int_{\theta=0}^{2 \pi} h_{1}(\overline{x} + \rho \cos \theta) \, h_{2}(\overline{x} + \rho \sin \theta)\, d \theta,\quad \overline{x} \in \mathbb{R}.
\end{equation*}
From \eqref{final:single} and the scaling, $S_{r}(f,g) = S_{1}(f \circ \delta_{r}, g \circ \delta_{r}) \circ \delta_{1/r}$, we obtain
\begin{align*}
& S_{r}(f,g)     (\underline{x}, \overline{x}) = S_{1} (f \circ \delta_{r}, g \circ \delta_{r}) \Big( \delta_{1/r} (\underline{x}, \overline{x})  \Big) \\
& \leq \int_{\rho=0}^{1}  {\scriptsize \rho \,  (  1 - \rho^2    )^{n-1}  }         \int_{ \theta=0}^{2 \pi} { \scriptsize M_{S}  (f \circ \delta_{r})  \big( r^{-1} \underline{x} , r^{-2} \overline{x} +  \rho \cos \theta \big)  }    { \scriptsize \  \Lambda  (g \circ \delta_{r})  \big( r^{-1}\underline{x}, r^{-2} \overline{x} +  \rho \sin \theta \big) } \, d \theta         \,             d \rho  \\
& = \int_{\rho=0}^{1}  {\scriptsize \rho \,  (  1 - \rho^2    )^{n-1}  }         \int_{ \theta=0}^{2 \pi} { \scriptsize M_{S}  (f \circ \delta_{r}) \circ \delta_{1/r}  \big(  \underline{x} ,  \overline{x} +  r^{2} \rho \cos \theta \big)  }    { \scriptsize \  \Lambda ( g \circ \delta_{r}) \circ \delta_{1/r}  \big( \underline{x},  \overline{x} + r^{2}  \rho \sin \theta \big) } \, d \theta         \,             d \rho  \\
 & = \int_{0}^{1}  {\scriptsize \rho \,  (  1 - \rho^2    )^{n-1}  }  B_{r^2 \rho}  \Big( { \scriptsize M_{S}  (f \circ \delta_{r}) \circ \delta_{1/r}  \big( \underline{x}, \cdot \big)  }  \, , \,   { \scriptsize \  \Lambda ( g  \circ \delta_{r}) \circ \delta_{1/r}) \big( \underline{x}, \cdot \big) } \Big)   (\overline{x})      \,             d \rho\\
  & \leq \int_{0}^{1}  {\scriptsize \rho \,  (  1 - \rho^2    )^{n-1}  }  B_{r^2 \rho} \Big( { \scriptsize M_{S}  f  \big( \underline{x}, \cdot \big)  }  \, , \,   { \scriptsize \  \Lambda  g   \big( \underline{x}, \cdot \big) } \Big)   (\overline{x})      \,             d \rho\\
  & \overset{ \rho \rightarrow r^{-2} \rho}{=} \frac{1}{r^{2}} \int_{0}^{r^{2}} {\textstyle \frac{\rho
  }{r^{2}} } \Big(1 - {\textstyle \big(\frac{\rho }{r^2} \big)^2 \Big) } ^{n-1}  B_{ \rho} \Big( { \scriptsize M_{S}  f  \big( \underline{x}, \cdot \big)  }  \, , \,   { \scriptsize \  \Lambda  g   \big( \underline{x}, \cdot \big) } \Big)  (\overline{x})\,d \rho.
 \end{align*}
As a consequence,
\begin{align} 
\mathcal{S}_{\rm{full}}(f, g)(\underline{x}, \overline{x})=\sup_{r>0} |S_{r}(f,g)     (\underline{x}, \overline{x})| &\lesssim \sup_{r>0} \frac{1}{r} \int_{0}^{r}   B_{ \rho}\Big( { \scriptsize M_{S}  f  \big( \underline{x}, \cdot \big)  }  \, , \,   { \scriptsize \  \Lambda  g   \big( \underline{x}, \cdot \big) } \Big)(\overline{x})\,\,d \rho\nonumber\\
&\lesssim \mathcal{B}(\scriptsize M_{S}  f  \big( \underline{x}, \cdot \big),   \scriptsize \Lambda  g   \big( \underline{x}, \cdot \big))(\overline{x}),  
\end{align} 
with the operator $\mathcal{B}$ acting only in the central direction, as usual. Interchanging the role of $f$ and $g$ we also obtain the estimate
\begin{align} 
\label{reverse}
\mathcal{S}_{\rm{full}}(f, g)(\underline{x}, \overline{x})\lesssim\mathcal{B}(\scriptsize   \Lambda f  \big( \underline{x}, \cdot \big),   \scriptsize  M_{S} g   \big( \underline{x}, \cdot \big))(\overline{x}). 
\end{align} 
This completes the proof. 
\end{proof}

\subsection{Proof of Theorem~\ref{full:main}}
Now we are at a position to prove our main result. 

\begin{proof}[Proof of Theorem~\ref{full:main}]
Let $\frac{2n}{2n-1} < p_{1} \leq \infty, 1< p_{2} \leq \infty$ satisfying $\frac{1}{p}=\frac{1}{p_1}+\frac{1}{p_2}.$ First using the boundedness of the operator $\mathcal{B}$ in the central variable from Theorem~\ref{euclidean:avg:thm}, followed by  H\"older's inequality we derive
\begin{align}\label{full: max 1}
   \nonumber\| \mathcal{S}_{\text{full}}(f, g) \|_{L^{p}(\Ha)}^{p} &\lesssim \int_{\C^n}\left( \int_{\R}|\mathcal{B}( M_{S}  f  \big( \underline{x}, \cdot \big),   \scriptsize \Lambda  g   \big( \underline{x}, \cdot \big))(\overline{x})|^{p}\,d\overline{x}\right)\,d\underline{x}\\ 
\nonumber &\lesssim \int_{\C^n} \left(\int_{\R} | M_{S}f(\underline{x}, \overline{x})|^{p_1}\, d\overline{x}\right)^{p/p_1} \left(\int_{\R} | \Lambda g(\underline{x}, \overline{x})|^{p_2}\, d\overline{x}\right)^{p/p_2} d\underline{x} \\
\nonumber &\lesssim \|M_{S}f\|_{L^{p_1}(\Ha)}^{p}\|\Lambda g\|_{L^{p_2}(\Ha)}^{p} \\
&\lesssim \| f \|_{L^{p_{1}}(\Ha)}^{p} \, \| g \|_{L^{p_{2}}(\Ha)}^{p},
\end{align}
since $M_{S}$ is bounded from $L^{p_{1}}(\Ha)$ to itself for $\frac{2n}{2n-1}<p_1\leq \infty,$ and $\Lambda$ maps $L^{p_{2}}(\Ha)$ to itself for $1<p_2 \leq \infty.$

Similarly, for $1< p_{1} \leq \infty, \frac{2n}{2n-1} < p_{2} \leq \infty,$ we can use the estimate \eqref{reverse}, to obtain
\begin{equation}\label{full: max 2}
   \| \mathcal{S}_{\text{full}}(f, g) \|_{L^{p}(\Ha)} \lesssim \| \Lambda f \|_{L^{p_{1}}(\Ha)} \| M_{S}g \|_{L^{p_{2}}(\Ha)} \,   \lesssim \| f \|_{L^{p_{1}}(\Ha)} \, \| g \|_{L^{p_{2}}(\Ha)}.
\end{equation}
Now a standard linearization technique and complex interpolation between \eqref{full: max 1} and \eqref{full: max 2} we obtain that (\ref{full: max 1}) holds when $(\frac{1}{p_{1}}, \frac{1}{p_{2}})$ belong to the open pentagon with the corners $(0,0)$, $(0,1)$, $(\frac{2n-1}{2n},1)$, $(1,\frac{2n-1}{2n})$ and $(1,0)$ with half open line segments $[O,A)$ and $[O,D)$ included.
    
\end{proof}

\section{Necessary conditions}
\label{sec:nece}
In this section, we construct examples to show some necessary conditions for the boundedness of $\mathcal{S}_{\rm{full}}.$ In fact, we prove a stronger statement. Let $\mathcal{S}_{\rm{loc}}$ be the local maximal function defined as
$$\mathcal{S}_{\rm{loc}}(f, g)(x)=\sup_{1\leq r\leq 2}|S_{r}(f, g)(x)|,$$
where $S_{r}$ is defined as in \eqref{def:single:scale}. Then we prove the following.

\begin{prop}
\label{nece:loc:new}
Let $n\geq 1.$ Let $1\leq p_1, p_2\leq \infty$ and $0<p<\infty.$ If we have \[\|\mathcal{S}_{\rm{loc}}(f, g)\|_{L^p(\Ha)}\leq C\|f\|_{L^{p_1}(\Ha)}\|g\|_{L^{p_2}(\Ha)}\] then we must have
\begin{align}
\label{Nece:2}
\frac{1}{p}\leq \frac{1}{p_1}+\frac{1}{p_2}\leq 2-\frac{1}{2n+1}=2-\frac{1}{d}.
\end{align}
\end{prop}

We break the proof into several lemmas. To begin with, let $\delta>0$ be a small positive number and $$\mathcal{E}:=\{x=(\underline{x}, \overline{x}): 1\leq |x|\leq 2,\quad \overline{x}>\kappa,\quad \|\underline{x}\|>\kappa \}.$$
Here, $\kappa>0$ is a fixed but very small number. For $x\in \mathcal{E},$ we denote $$r=r(x):=(4\|\underline{x}\|^4+2|\overline{x}|^2)^{1/4}, \, \text{then}\,\, r(x)\sim 1.$$ 

\begin{lemma}
Fix $x\in \mathcal{E}.$ Set $$E^{\delta}_{x, 1}:=\{\rho\in [0, 1]: \big|\frac{\sqrt{2}\overline{x}}{r(x)^2}-\rho\big|\lesssim \delta\}.$$ Then we have the estimate 
\begin{align}
\label{imp:nece:1}
\big|r(x)(1-\rho^2)^{\frac{1}{4}}-\sqrt{2}\|\underline{x}\|\big|\lesssim \delta. 
\end{align}
\end{lemma}
\begin{proof}
Observe that $$r(x)^4(1-\rho^2)-4\|\underline{x}\|^4=r(x)^4-4\|\underline{x}\|^4-r(x)^4\rho^2=2|\overline{x}|^2-r(x)^4\rho^2=(\sqrt{2} |\overline{x}|-r(x)^2\rho)(\sqrt{2} |\overline{x}|+r(x)^2\rho).$$ 
As $\rho\in E_{x, 1},$ $|r(x)^4(1-\rho^2)-4\|\underline{x}\|^4|\lesssim |\sqrt{2} |\overline{x}|-r(x)^2\rho|\lesssim \delta.$ The last estimate further implies that 
\begin{align*}
\big|r(x)(1-\rho^2)^{\frac{1}{4}}-\sqrt{2}\|\underline{x}\|\big|\lesssim \delta. 
\end{align*}
\end{proof}

\begin{lemma}
\label{z}
Fix $x\in \mathcal{E}$ and denote $\widehat{\underline{x}}:=\underline{x}/ \| \underline{x} \|.$ Set $E^{\delta}_{x, 2}:=\{z_1\in B^{2n}(0, 1): \big|\widehat{\underline{x}}-\sqrt{2}z_{1}\big|\lesssim \delta\}.$ Then we have the estimate \begin{align}
\bigg|\sqrt{1-\|z_1\|^2}-\frac{1}{\sqrt{2}}\bigg|\lesssim \delta.    \label{imp:nece:4}
\end{align}   
\end{lemma}

\begin{proof}
For $z_{1}\in E^{\delta}_{x, 2}$, by triangle inequality, we obtain that $|\frac{1}{\sqrt{2}}-\|z_{1}\||\lesssim \|\hat{\underline{x}}-\sqrt{2} z_1\|=O(\delta),$ and this in turn yields
\begin{align}
\nonumber&\big|\sqrt{1-\|z_1\|^2}-\frac{1}{\sqrt{2}}\big| \big|\sqrt{1-\|z_1\|^2}+\frac{1}{\sqrt{2}}\big|\\
&=\big(1-\|z_{1}\|^2-\frac{1}{2}\big)=\big(\frac{1}{2}-\|z_{1}\|^2\big)=\big|\frac{1}{\sqrt{2}}-\|z_1\|\big| \big|\frac{1}{\sqrt{2}}+\|z_1\|\big|=O(\delta).
\label{ex:sharp:1}
\end{align}
Consequently, \begin{align*}
\bigg|\sqrt{1-\|z_1\|^2}-\frac{1}{\sqrt{2}}\bigg|\lesssim \delta.
\end{align*}
\end{proof}

Henceforth, $B^{e}(0, c\,\delta)$ denotes the Euclidean ball centered at the origin with radius $c\delta$ in $\R^{2n+1},$ and $c>0$ is a fixed constant, chosen to satisfy the upcoming requirements of the proofs.

\begin{lemma}
\label{f}
Fix $\delta>0,$ and $f=\chi_{B^{e}(0, c \delta)}.$ Let $x\in \mathcal{E}.$ Moreover, let $\rho\in E^{\delta}_{x, 1},$ $z_{1}\in E^{\delta}_{x, 2},$ and
$$\theta\in F^{\delta}_{x}:=\{\theta\in [0, 2\pi]: \big|\big(\frac{\overline{x}}{\sqrt{2}|\overline{x}|}, \frac{\overline{x}}{\sqrt{2}|\overline{x}|}\big)-\big(\cos \theta, \sin \theta\big) \big|\lesssim \delta\}.$$
Then \begin{align}
\label{1st:comp}
 f(x . ( r(x)(1- \rho^2)^{1/4}z_{1}, r(x)^2\rho \cos \theta)^{-1})=1.   
\end{align}
\end{lemma}

\begin{proof}
It is enough to show $x . ( r(x)(1- \rho^2)^{1/4}z_{1}, r(x)^2\rho \cos \theta)^{-1}\in B^{e}(0, c\, \delta).$ We denote the coordinates of $x . ( r(x)(1- \rho^2)^{1/4}z_{1}, r(x)^2\rho \cos \theta)^{-1}$ as $(\underline{y}, \bar{y}),$ that is, 
\begin{align*}
 &\underline{y}=\underline{x}-r(x)(1-\rho^2)^{1/4}z_{1},\\
 &\bar{y}=\overline{x}-r(x)^2\rho \cos\theta-\frac{r(x)(1-\rho^2)^{1/4}}{2}\Im(\underline{x}.\bar{z_1}).
\end{align*}

Now, observe that
\begin{align}
\label{imp:nece:2}
\|\underline{y}\|\leq \|\underline{x}-\sqrt{2} \|\underline{x}\|z_{1}\|+ \|\sqrt{2} \|\underline{x}\|z_{1}- r(x)(1-\rho^2)^{1/4}z_{1}\|\lesssim \delta    
\end{align}
as $z_1\in E^{\delta}_{x, 2}$ and \eqref{imp:nece:1}. The estimate \eqref{imp:nece:2} further implies that 
$$|\Im (\underline{x}\bar{z_{1}})|\lesssim \|\underline{y}\|\lesssim \delta.$$
Therefore, in order to prove $|\bar{y}|\lesssim \delta,$ we need only to show that 
\begin{align}
\label{imp:nece:3}
\big|\overline{x}-r(x)^2\rho \cos\theta\big|\lesssim \delta.    
\end{align}
To prove \eqref{imp:nece:3}, observe that
\begin{align}
\nonumber &\big|\overline{x}-r(x)^2\rho \cos\theta\big|\\
\nonumber&\lesssim \big|\overline{x}-\sqrt{2}|\overline{x}| \cos \theta\big|+\big|\sqrt{2}|\overline{x}| \cos \theta-r(x)^2\rho \cos\theta\big|\\
\label{angle}&\lesssim \big|\overline{x}-\sqrt{2}|\overline{x}| \cos \theta\big|+\big|r(x)^2\rho -\sqrt{2} |\overline{x}|\big|\lesssim \delta,
\end{align}
as $\rho\in E^{\delta}_{x, 1}$ and $\theta\in F_{x}^{\delta}.$ Hence $(\underline{y}, \bar{y})\in B^{e}(0, c\delta)$ and thus $f(\underline{y}, \bar{y})=1.$
\end{proof}

The final lemma we require reads as follows. For $x\in \mathcal{E},$ set  $$E^{\delta}_{x, 3}:=\{z_{2}\in S^{2n-1}: \big|\widehat{\underline{x}}-z_{2}\big|\lesssim \delta\}.$$

\begin{lemma}
\label{g}
Fix $\delta>0,$ and $g=\chi_{B^{e}(0, c\delta)}.$ Let $x\in \mathcal{E}.$ Moreover, let $\rho\in E^{\delta}_{x, 1},$ $z_{1}\in E^{\delta}_{x, 2},$ and
$$\theta\in F^{\delta}_{x}:=\{\theta\in [0, 2\pi]: \big|\big(\frac{\overline{x}}{\sqrt{2}|\overline{x}|}, \frac{\overline{x}}{\sqrt{2}|\overline{x}|}\big)-\big(\cos \theta, \sin \theta\big) \big|\lesssim \delta\},\quad z_{2}\in E^{\delta}_{x, 3}.$$
Then \begin{align}
\label{2st:comp}
g ( x . (r(x)(1- \rho^2)^{1/4} {\scriptsize \sqrt{1-\|z_{1}\|^2}} z_{2}, r(x)^2\rho \sin \theta  )^{-1} )=1.
\end{align}
\end{lemma}

\begin{proof}
We denote $(\underline{a}, \overline{a})=( x . (r(x)(1- \rho^2)^{1/4} {\scriptsize \sqrt{1-\|z_{1}\|^2}} z_{2}, r(x)^2\rho \sin \theta  )^{-1} ).$ Then
 \begin{align*}
&\underline{a}:= \underline{x}-r(x)(1- \rho^2)^{1/4} {\scriptsize \sqrt{1-\|z_{1}\|^2}} z_{2}\\
&\bar{a}:=\overline{x}-r(x)^2\rho \sin \theta-\frac{r(x)(1- \rho^2)^{1/4} {\scriptsize \sqrt{1-\|z_{1}\|^2}}}{2}\Im(\underline{x}\bar{z_{2}}).
\end{align*}

By the triangle inequality, we observe
\begin{align}
\nonumber&\|\underline{a}\|=\|\underline{x}-r(x)(1- \rho^2)^{1/4} {\scriptsize \sqrt{1-\|z_{1}\|^2}} z_{2}\|\\
\nonumber &\lesssim \|\underline{x}-\sqrt{2} \|\underline{x}\| {\scriptsize \sqrt{1-\|z_{1}\|^2}} z_{2}\|+\|\sqrt{2} \|\underline{x}\| {\scriptsize \sqrt{1-\|z_{1}\|^2}} z_{2}-r(x)(1- \rho^2)^{1/4} {\scriptsize \sqrt{1-\|z_{1}\|^2}} z_{2}\|\\
\nonumber &\lesssim \|\underline{x}-\sqrt{2} \|\underline{x}\| {\scriptsize \sqrt{1-\|z_{1}\|^2}} z_{2}\|+\|\sqrt{2} \|\underline{x}\|-r(x)(1- \rho^2)^{1/4}\|\\
\nonumber &\lesssim \|\underline{x}-\|\underline{x}\| z_{2}\|+\|\|\underline{x}\| z_{2}-\sqrt{2} \|\underline{x}\| {\scriptsize \sqrt{1-\|z_{1}\|^2}} z_{2}\|+\|\sqrt{2} \|\underline{x}\|-r(x)(1- \rho^2)^{1/4}\|\\
&\lesssim \|\widehat{\underline{x}}- z_{2}\|+\bigg|\frac{1}{\sqrt{2}}-{\scriptsize \sqrt{1-\|z_{1}\|^2}}\bigg|+\|\sqrt{2} \|\underline{x}\|-r(x)(1- \rho^2)^{1/4}\|\lesssim \delta,\label{2nd:comp:final}
\end{align}
where the last inequality follows as $z_{2}\in E^{\delta}_{x, 3},$ \eqref{imp:nece:1} and \eqref{imp:nece:4}. The estimate \eqref{2nd:comp:final} also implies 
$$|\Im (\underline{x} . \bar{z_2} )|\lesssim \|\underline{a}\|\lesssim \delta.$$

To show that $|\overline{a}|\lesssim \delta,$ we need to show $|\overline{x}-r(x)^2\rho \sin \theta|\lesssim \delta,$ but this follows along the same lines as in \eqref{angle}. This concludes $g(\underline{a}, \overline{a})=1.$
\end{proof}

Now we are in a position to prove the Proposition~\ref{nece:loc:new}.
 
\begin{proof}[Proof of Proposition~\ref{nece:loc:new}]
The condition  $\frac{1}{p} \leq \frac{1}{p_1}+ \frac{1}{p_2}$ is even necessary for the boundedness of $S_{r}$ since each $S_{r}$ is a local operator satisfying $\tau_{y} (S_{r}(f,g) )=S_{r}(\tau_{y}f, \tau_{y}g),$ where $\tau_{y}$ is a left translation defined as $\tau_{y}(h)(x):= h(y.x)$ for all $x \in \Ha.$

Now we prove the second inequality in \eqref{Nece:2}. Let $\delta>0$ be any small positive number and we choose the functions \begin{align*}
   f=\chi_{B^{e}(0,\, c\,\delta)},\, \text{and}\quad g=\chi_{B^{e}(0,\, c\,\delta)}.  
\end{align*}

The proof will follow from the estimate
\begin{align}
\label{main:nece}
\mathcal{S}_{\rm{loc}}(f, g)(x)\gtrsim \delta^{4n+1}
\end{align}
for $x\in \mathcal{E}.$ Since $|\mathcal{E}|\sim 1,$ once we prove \eqref{main:nece}, the condition \eqref{Nece:2} follows since
$$\delta^{4n+1}|\mathcal{E}|\lesssim \|\mathcal{S}_{\rm{loc}}\|_{L^p(\Ha)}\lesssim \|f\|_{L^{p_1}(\Ha)} \|g\|_{L^{p_2}(\Ha)}\lesssim \delta^{\frac{2n+1}{p_1}}\, \delta^{\frac{2n+1}{p_2}}.$$
Letting $\delta\to 0,$ it follows that $$\frac{1}{p_1}+\frac{1}{p_2}\leq \frac{4n+1}{2n+1}=2-\frac{1}{2n+1}.$$

Now we proceed to prove \eqref{main:nece}. Fix some $x=(\underline{x}, \overline{x})\in \mathcal{E}.$ Recall the expression derived in \eqref{slice:final:1} stating
\begin{align}
\nonumber S(f,g)(\underline{x}, \overline{x})& \sim \int_{\rho=0}^{1} \int_{ \theta=0}^{2 \pi} {\scriptsize \rho \,  (  1 - \rho^2    )^{n-1}  } \int_{z_{1} \in B^{2n}(0,1)}   {\scriptsize  f ( x . ( (1- \rho^2)^{1/4}z_{1}, \rho \cos \theta)^{-1} ) } \\
& \ \ \ \ \int_{z_{2} \in S^{2n-1}} {\scriptsize g ( x . ( (1- \rho^2)^{1/4} {\scriptsize \sqrt{1-\|z_{1}\|^2}} z_{2}, \rho \sin \theta  )^{-1} ) }   d \sigma_{2n-1}(z_{2}) \, d \rho \,  (1-\|z_{1}\|^2)^{n-1} d z_{1}\, d\theta\label{slice:nece}.
\end{align}

Using the scaling property $S_{r}(f,g) = S(f \circ \delta_{r}, g \circ \delta_{r}) \circ \delta_{1/r},$ \eqref{slice:nece} implies that
\begin{align}
\nonumber S_{r}(f,g)(\underline{x}, \overline{x})& \sim \int_{ \theta=0}^{2 \pi}\int_{\rho=0}^{1}  {\scriptsize \rho \,  (  1 - \rho^2    )^{n-1}  } \int_{z_{1} \in B^{2n}(0,1)}   {\scriptsize  f ( x . ( r(1- \rho^2)^{1/4}z_{1}, r^2\rho \cos \theta)^{-1} ) } \\
& \ \ \  \int_{z_{2} \in S^{2n-1}} {\scriptsize g ( x . (r(1- \rho^2)^{1/4} {\scriptsize \sqrt{1-\|z_{1}\|^2}} z_{2}, r^2\rho \sin \theta  )^{-1} ) }   d \sigma_{2n-1}(z_{2}) \, d \rho \,  (1-\|z_{1}\|^2)^{n-1} d z_{1}d\theta\label{slice:nece:2}.
\end{align}
We aim to prove $$S_{r(x)}(f, g)(x)\gtrsim \delta^{4n+1},\quad \text{where}\quad r(x):=(4\|\underline{x}\|^4+2|\overline{x}|^2)^{1/4}.$$

For $x\in \mathcal{E},$ recall the sets 
\begin{align*}
&E^{\delta}_{x, 1}:=\{\rho\in [0, 1]: \big|\frac{\sqrt{2}\overline{x}}{r(x)^2}-\rho\big|\lesssim \delta\}\\
&E^{\delta}_{x, 2}:=\{z_1\in B^{2n}(0, 1): \big|\widehat{\underline{x}}-\sqrt{2}z_{1}\big|\lesssim \delta\}\\
&E^{\delta}_{x, 3}:=\{z_{2}\in S^{2n-1}: \big|\widehat{\underline{x}}-z_{2}\big|\lesssim \delta\}\\
&F^{\delta}_{x}:=\{\theta\in [0, 2\pi]: \big|\big(\frac{\overline{x}}{\sqrt{2}|\overline{x}|}, \frac{\overline{x}}{\sqrt{2}|\overline{x}|}\big)-\big(\cos \theta, \sin \theta\big) \big|\lesssim \delta\}.
\end{align*}

Then
\begin{align}
\nonumber& S_{r(x)}(f,g)(\underline{x}, \overline{x})\\
\nonumber& \gtrsim  \int_{\theta\in F_{x}^{\delta}} \int_{\rho\in E^{\delta}_{x, 1}}  {\scriptsize \rho \,  (  1 - \rho^2    )^{n-1}  } \int_{z_{1}\in E^{\delta}_{x, 2}}   {\scriptsize  f ( x . ( r(1- \rho^2)^{1/4}z_{1}, r^2\rho \cos \theta)^{-1} )} \\
\nonumber& \ \ \ \ \ \ \ \ \ \ \ \ \int_{z_{2}\in E_{x, 3}^{\delta}} {\scriptsize g ( x . (r(1- \rho^2)^{1/4} {\scriptsize \sqrt{1-\|z_{1}\|^2}} z_{2}, r^2\rho \sin \theta  )^{-1} ) }   d \sigma_{2n-1}(z_{2}) \, d \rho \,  (1-\|z_{1}\|^2)^{n-1} d z_{1}\, d\theta.
\end{align}
By Lemma~\ref{f} and Lemma~\ref{g}, we have $f ( x . ( r(1- \rho^2)^{1/4}z_{1}, r^2\rho \cos \theta)^{-1} )=1,$ and $ g ( x . (r(1- \rho^2)^{1/4} {\scriptsize \sqrt{1-\|z_{1}\|^2}} z_{2}, r^2\rho \sin \theta  )^{-1} )=1,$ respectively. Also, $(1-\|z_1\|^2)\sim 1$ by Lemma~\ref{z}. Since $x\in \mathcal{E},$ and $\rho\in E^{\delta}_{x, 1},$ choosing $\delta$ sufficiently small we can conclude that $\rho\geq c_1>0,$ for some small but fixed constant $c_1.$ Moreover, the inequality \eqref{imp:nece:1} implies that $(1-\rho^2)\geq c_2,$ for some small but fixed constant $c_2.$ Combining these estimate we derive
\begin{align*}
S_{r(x)}(f,g)(\underline{x}, \overline{x})&\gtrsim |F_{x}^{\delta}| |E_{x, 1}^{\delta}| |E_{x, 2}| \sigma_{2n-1}(E_{x, 3})\gtrsim {\delta} \delta \delta^{2n} \delta^{2n-1}\simeq \delta^{4n+1}.     
\end{align*}
Here, of course, $|F_{x}^{\delta}|$ is the arc-length measure of $F_{x}^{\delta}\subset S^{1}$ and $|F_{x}^{\delta}|\sim \delta,$ $|E_{x, 1}^{\delta}|$ is the Lebesgue measure of $E_{x, 1}^{\delta}\subset \R$ and $|E_{x, 1}^{\delta}|\sim \delta,$ $|E_{x, 2}|$ is the Lebesgue measure $E_{x, 2}$ on $\R^{2n}$ with $|E_{x, 2}|\sim \delta^{2n},$ and finally $\sigma_{2n-1}(E_{x, 3})\sim \delta^{2n-1}.$ This completes the claim \eqref{main:nece} and hence completes the proof of Proposition~\ref{nece:loc:new}.

\end{proof}

\section{Multilinear generalizations}
\label{multi}
Our arguments can also be generalized to the $k$-linear setting for $k\geq 2.$ We briefly outline it here. Let $\mf{G}_k \subset (\R^{2n+1})^k$ be the surface defined as
\begin{align*}
\{ (z_1, t_1, \dots, z_k, t_k): \Phi(z_1, t_1, \dots, z_k, t_k) := \big( \sum_{j=1}^k \|z_j\|^2 \big)^2 + \sum_{j=1}^k t_j^2 - 1 = 0\},
\end{align*}
where $\vec{\mathfrak{z}} = (z_1, \dots, z_k) \in (\mathbb{R}^{2n})^k,$ and $\vec{\mathfrak{t}} =(t_1, \dots, t_k) \in \mathbb{R}^k$ is the central variable. Also, we denote $\vec{\mf{z}'} = (z_2, \dots, z_k).$ Define the $k$-linear Kor\'anyi maximal function by
\begin{align}
 &\mathcal{S}_{\rm{full}}^{k}(f_1, \cdots, f_{k}):=\sup_{r>0}|S_{r}^{k}(f_1, \cdots, f_{k})(x)|,\\
& \text{where}\quad S_{r}^{k}(f_1, \cdots, f_{k})(x)=\int_{\mf{G}_{k}} \,\prod_{i=1}^{k} f_{i} \big( x \cdot \delta_{r}(z_{i},t_{i}) \big)\,d \sigma(\vec{\mf{z}}, \vec{\mf{t}}),\quad x\in \Ha,
\end{align}
and $d\sigma$ is the normalized induced Lebesgue measure on $\mathfrak{G}_{k}$ in $(\R^{2n+1})^k.$ Using \eqref{slicing:formula} and arguing exactly as in the proof of Proposition~\ref{slice:eqn}, we derive
\begin{align}
\label{multi:step:1}
S_r^{k}(f_1, \dots, f_k)(x) \sim \int_{\|\vec{\mf{t}}\|<1} (1 - \|\vec{\mf{t}}\|^2)^{\frac{kn-2}{2}} \int_{\vec{\mf{z}} \in S^{2kn-1}} \prod_{j=1}^k f_j\left(x \cdot \left(r (1-\|\vec{\mf{t}}\|^2)^{1/4} z_j, r^2 t_j\right)\right) d\sigma_{2kn-1}(\vec{\mf{z}}) \, d\vec{\mf{t}}.
\end{align}
A repeated application of the formula 
\begin{align*}
\int_{S^{2kn-1}} G \, d\sigma_{2kn-1} = \int_{z_1 \in B^{2n}(0,1)} \left[ \int_{S^{2(k-1)n-1}} G\left(z_1, \sqrt{1-\|z_1\|^2}\vec{\mf{z'}}\right) d\sigma_{2(k-1)n-1}(\vec{\mf{z'}}) \right] (1-\|z_1\|^2)^{(k-1)n-1} \, dz_1,
\end{align*}
will lead us to the following pointwise estimate
\begin{align}
\label{multi:point}
|S_r^{k}(f_1, \dots, f_k)(\underline{x}, \overline{x})| \lesssim \int_0^1 \varphi(\rho)\,{B}_{r^2\rho}^{k} \Big( M_S f_1(\underline{x}, \cdot), \Lambda f_2(\underline{x}, \cdot), \dots, \Lambda f_k(\underline{x}, \cdot) \Big)(\overline{x}) \, d\rho,
\end{align}
where $\varphi$ is a polynomial in $\rho,$ and ${B}^{k}_{\rho}$ is the $k-$linear averaging operator acting on the central variable $\overline{x},$ defined by
$${B}^{k}_{\rho}(g_1,\cdots, g_{k})(\overline{x})=\int_{S^{k-1}}\prod_{i=1}^{k}g_{i}( \overline{x}-\rho y_{i})\, d\sigma_{k-1}(y_1, \cdots, y_k)\,,\quad \overline{x}\in \R,$$
where $g_i \in C(\R).$  A further change of variables $\rho\to r^{-2}\rho$ in \eqref{multi:point} yields
\begin{align}
\label{multi:point:1}
\mathcal{S}^{k}_{\rm{full}}(f_1, \cdots, f_{k})(\underline{x}, \overline{x}) \lesssim \mathcal{B}^{k}\Big( M_S f_1(\underline{x}, \cdot), \Lambda f_2(\underline{x}, \cdot), \dots, \Lambda f_k(\underline{x}, \cdot) \Big)(\overline{x}),
\end{align}
where 
\begin{align}
\mathcal{B}^{k}(g_{1}, \cdots, g_{k})(\overline{x})= \sup_{r>0}\big|\frac{1}{r}\int_{0}^{r} {B}^{k}_{\rho}(g_1,\cdots, g_{k})(\overline{x})\,\, d\rho\big|,\quad \overline{x}\in \R.   
\end{align}

A few remarks on the averages $B_{\rho}^{k}$ are in order. These operators were introduced by D. Oberlin in \cite{Obe88}, where $L^{p}(\mathbb{R}) \times \cdots \times L^{p}(\mathbb{R}) \to L^{q}(\mathbb{R})$ boundedness results were established, and later sharpened by Bak and Shim in \cite{Bak-Sim}. More recently, the authors of \cite{SS21} studied the boundedness of $B_{\rho}^{k}$ for general exponents $1 \leq p_1,\ldots,p_k,r \leq \infty$, such that the boundedness $B_{\rho}^{k} : L^{p_1}(\mathbb{R}) \times \cdots \times L^{p_k}(\mathbb{R}) \to L^{r}(\mathbb{R})$ holds. However, for the maximal operator $\mathcal{B}^{k},$ following our arguments in Theorem~\ref{euclidean:avg:thm} we can derive the following result.

\begin{theorem}
\label{euclidean:avg:thm:multi}
Let $1<p_1, \cdots, p_k\leq \infty$ and $\frac{1}{p}=\frac{1}{p_1}+\cdots+\frac{1}{p_k}.$ Then $\mathcal{B}^k$ maps $L^{p_1}(\mathbb{R})\times \cdots  L^{p_k}(\mathbb{R})$ boundedly to $L^{p}(\mathbb{R}).$ If $p_i=1$ for some $1\leq i\leq k$, the corresponding weak-type estimates hold.    
\end{theorem}

Arguing exactly as in \eqref{multi:point} one can also obtain
\begin{align}
\label{multi:point:gen}
\mathcal{S}^{k}_{\rm{full}}(f_1, \cdots, f_{k})(\underline{x}, \overline{x}) \lesssim \mathcal{B}^{k}\Big( \Lambda f_1(\underline{x}, \cdot),\cdots, \Lambda f_{i-1}(\underline{x}, \cdot),\, M_{S}f_{i}(\underline{x}, \cdot), \Lambda f_{i+1}(\underline{x}, \cdot), \cdots,\Lambda f_k(\underline{x}, \cdot) \Big)(\overline{x}),\quad x=(\underline{x}, \overline{x})\in \Ha,
\end{align}
for any $i=1,\cdots, k.$ Using Theorem~\ref{euclidean:avg:thm:multi}, \eqref{multi:point:gen}, and the boundedness of $M_{S},$ and $\Lambda$ we conclude with the following multilinear generalization of Theorem~\ref{full:main}. 

\begin{theorem}
\label{full:main:multi}
Let $n\geq 2.$ Suppose $1< p_1, p_2, \cdots, p_{k}\leq \infty,$  $\frac{1}{p}=\frac{1}{p_1}+
\frac{1}{p_2}+\cdots+\frac{1}{p_k}$ and $\frac{1}{p}<\frac{2nk-1}{2n}.$ Then $\mathcal{S}_{\mathrm{full}}^{k}$ is bounded from $\prod_{i=1}^{k}L^{p_{i}}(\Ha)$ to $L^p(\Ha).$
\end{theorem}

\medskip

\subsection*{Acknowledgements} AG gratefully acknowledges the support by the Industrial Consultancy and Sponsored Research (IC \& SR), Indian Institute of Technology Madras for the New Faculty Initiation Grant RF25261459MANFIG009296.

\subsection*{Author Contributions} All authors contributed equally. 

\subsection*{Data availability} Data sharing is not applicable to this article as no new data were created or analyzed in this study.

\section*{Declarations}

\subsection*{Tool and computational resource disclosure} No AI tools were used in the production of this paper.

\subsection*{Conflict of interest:} The authors declare that there is no conflict of interest.

\subsection*{Competing interests} The authors declare no competing interests.

\end{document}